\documentclass[11pt]{article}
\usepackage{graphicx}
\newtheorem{theorem}{Theorem}
\newtheorem{lemma}[theorem]{Lemma}
\newtheorem{sublemma}[theorem]{Sublemma}
\newtheorem{proposition}[theorem]{Proposition}
\newtheorem{conjecture}[theorem]{Conjecture}
\newtheorem{corollary}[theorem]{Corollary}

\newtheorem{definition}[theorem]{Definition}

\newtheorem{notation}[theorem]{Notation}

\begin{document}

\title{A result on resolutions of Veronese embeddings}
\author{  Elena Rubei }
\date{\hspace*{2cm}}
\maketitle

\vspace{-1cm}

{\small
 {\bf Abstract.} 
This paper deals with syzygies  of the ideals of the Veronese embeddings.
By Green's Theorem we know  that  ${\cal O}_{{\bf P}^n}(d)$   
satisfies Green-Lazarsfeld's Property $N_p$ 
 $\forall d \geq p$,  $\forall n $. 
By Ottaviani-Paoletti's theorem  if $n \geq 2$,
  $d \geq 3$ and $ 3d -2 \leq p$ then
${\cal O}_{ {\bf P}^{n} } (d)$ does not  satisfy Property $N_{p} $.
The cases $n \geq 3$, $d \geq 3$, $d < p < 3d -2$ are still open (except
 $n=d=3$).

Here we deal with one of these cases, namely
 we prove that  ${\cal O}_{{\bf P}^n} (3)$
 satisfies Property $N_4$ $\forall n$. 

Besides we prove that   
 ${\cal O}_{{\bf P}^n}(d)$   satisfies $N_p$  $\forall n \geq  p$ 
iff  ${\cal O}_{{\bf P}^{p}}(d)$  satisfies $N_p$.   }

\section{Introduction}

Let  $ L $ be a very ample  line bundle on a smooth
complex projective variety $Y$ and let  $\varphi_{L}: Y \rightarrow
{\bf P}(H^{0}(Y, L)^{\ast})$ be  the map  associated to $L$.
We  recall the definition of Property $N_{p}$ of Green-Lazarsfeld,
 studied for the first time  by Green in  \cite{Green1}
(see also \cite{G-L}, \cite{Green2}):

{\em
let $Y$ be a smooth complex projective variety and let $L$ be a very ample
 line bundle on $Y$ defining an embedding $\varphi_{L}: Y \hookrightarrow
{\bf P}={\bf P}(H^{0}(Y,L)^{\ast })$;
set $S= S(L)= \oplus_n Sym^{n} H^{0}(L) $,
the homogeneous coordinate ring of the projective space ${\bf P}$, and
consider
the graded $S$-module  $G=G(L)= \oplus_{n} H^{0}(Y, L^{n})$; let $E_{\ast}$
\[ 0 \longrightarrow E_{l} \longrightarrow E_{l-1} \longrightarrow 
... \longrightarrow E_{0} \longrightarrow G \longrightarrow 0\] 
be a minimal graded free resolution of $G$;
the line bundle $L$ satisfies Property
$ N_{p}$ ($p \in {\bf N}$) iff

\hspace{1cm} $E_{0}= S$

\hspace{1cm} $E_{i}= \oplus S(-i-1)$    \hspace{1cm}
for $1 \leq i \leq p $.}

(Thus $L$ satisfies Property $N_{0}$  iff  $Y \subset
{\bf P}(H^{0}(L)^{\ast })$ is projectively normal;
 $L$ satisfies Property $N_{1}$
iff $L$ satisfies $N_{0}$
and the homogeneous ideal $I$ of $Y \subset
{\bf P}(H^{0}(L)^{\ast })$ is generated by quadrics;
$L$ satisfies  $N_{2}$  iff
 $L$ satisfies  $N_{1}$ and the module of syzygies among
quadratic
generators $Q_{i} \in I$ is spanned by relations of the form
 $\sum L_{i}Q_{i}=0$, where $L_{i}$  are linear polynomials;
and so on.)

In this paper we will consider the case of Veronese embedding i.e. 
 the case  $Y= {\bf P}^n$,
 $L ={\cal O}(d)$. 
Among the papers on syzygies 
in this case we quote  \cite{B-M}, \cite{Green1}, \cite{O-P},
\cite{J-P-W}, \cite{Las}, \cite{P-W}. 
Two of the  most important results are:

\begin{theorem} \label{Gr} ({\bf Green}) \cite{Green1}.
Let $d,p \in {\bf N} $. If $d \geq p$
then ${\cal O}_{{\bf P}^{n}} (d)$ 
  satisfies Property $N_p$.
\end{theorem}

\begin{theorem} \label{OP}  ({\bf Ottaviani-Paoletti}) \cite{O-P}.
If $n \geq 2$,
  $d \geq 3$ and $ 3d -2 \leq p$ then
${\cal O}_{ {\bf P}^{n} } (d)$ does not  satisfy Property $N_{p} $.
\end{theorem}

\smallskip

\hspace{5cm}
\includegraphics[scale=0.24]{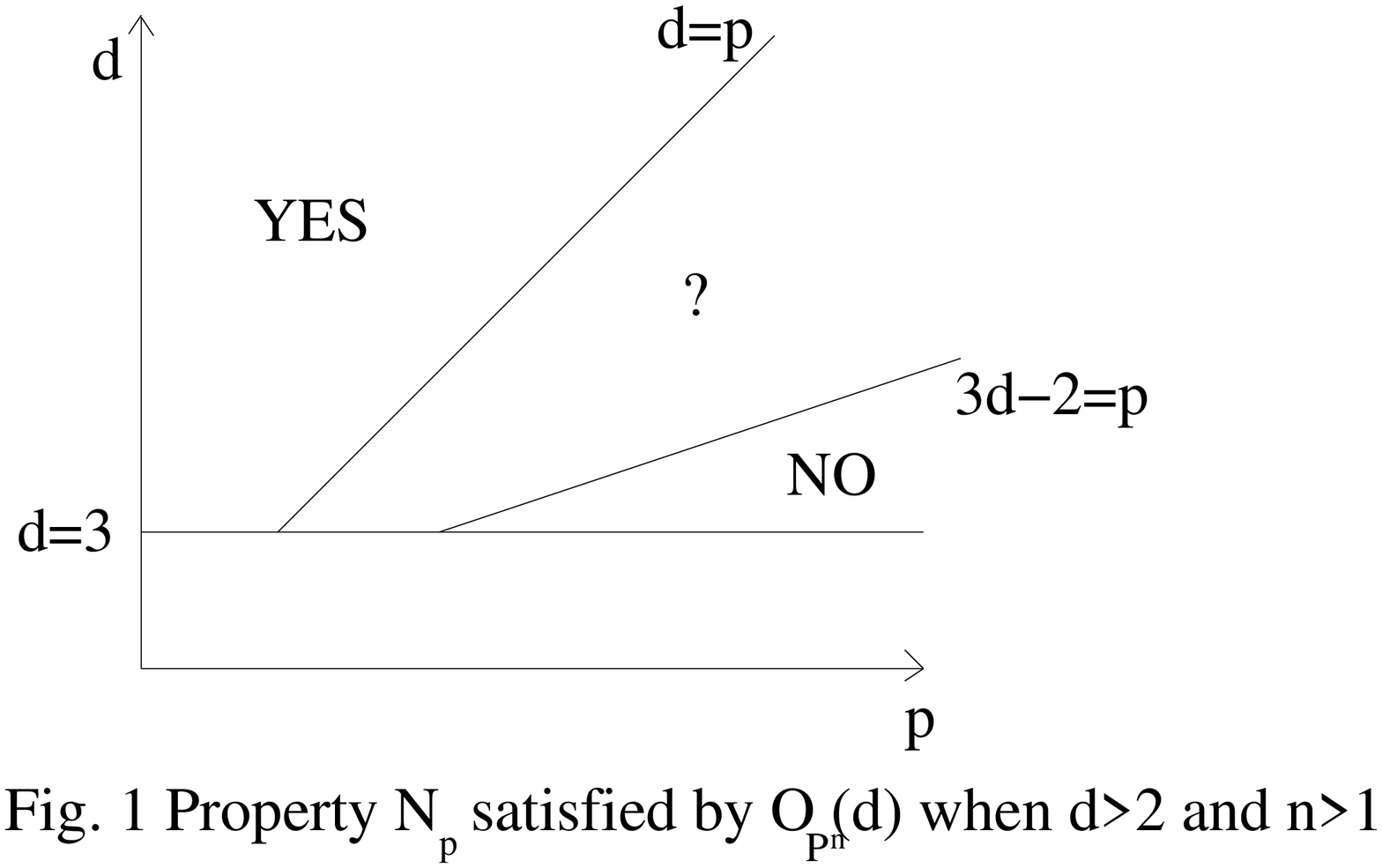}

\smallskip

In \cite{O-P} Ottaviani and Paoletti conjectured:

\begin{conjecture} \label{cOP}  ({\bf Ottaviani-Paoletti}) \cite{O-P}.
Let $n \geq 2$,  $d \geq 3$. The bundle
${\cal O}_{ {\bf P}^{n} } (d)$ satisfies Property $N_{p} $ iff $
p < 3d -2$.
\end{conjecture}

The following theorems and  proposition   show that in the case $n=2$
and in the case $n=d=3$  Conjecture \ref{cOP} is true:

\begin{theorem} \label{JPW} ({\bf Josefiak-Pragacz-Weyman}) \cite{J-P-W}.
Let $n \geq 3$. Then  ${\cal O}_{ {\bf P}^{n} } (2)$ 
satisfies  $N_{p} $ iff $p \leq 5$.
\end{theorem}

\begin{theorem} \label{G-B} ({\bf Green-Birkenhake}) \cite{Green1}, \cite{Bir}.
The bundle
${\cal O}_{ {\bf P}^{2} } (d)$ satisfies Property $N_{3d-3} $.
\end{theorem}

(It is well known that  ${\cal O}_{ {\bf P}^{2} } (2)$ 
satisfies Property $N_p $ $ \forall p$, see for instance \cite{O-P}.)

\begin{proposition} \label{OP3}  ({\bf Ottaviani-Paoletti}) \cite{O-P}.
The bundle ${\cal O}_{ {\bf P}^{3} } (3)$ satisfies Property $N_{6} $.
\end{proposition}

We recall also that the following Ein-Lazarsfeld's result implies 
Green's Theorem:

\begin{theorem}[Ein-Lazarsfeld]\label{t4}  \cite{E-L}
Let $Y$ be  a smooth complex projective
variety of dimension $n$; let $ A $ be a very ample line bundle on $Y$
and $ B$  a numerically effective line bundle on $Y$; then 
$ K_{Y} \otimes A^{n+1+ p } \otimes  B $ satisfies Property  $N_{p}$.          
Besides, if $ (Y, A, B) \neq ({\bf P}^{n}, {\mathcal O}_{{\bf P}^{n}}(1), 
{\mathcal O}_{{\bf P}^{n}})$ and $ p \geq 1$
 then $ K_{Y} \otimes A^{n + p} \otimes B $ satisfies Property $N_{p}$.
\end{theorem}

Also the following result implies Green's Theorem (taking $Y = {\bf P}^n $ 
and  $M= {\cal O}(1)$).
 
\begin{theorem} [Rubei]\label{1} \cite{Ru1}
 Let $Y$ be a smooth complex projective variety and $M$  a line bundle 
on $Y$.
 If $M$ satisfies Property $ N_{p}$ and $ d \geq p $ 
then $M^{d}$ satisfies Property
$N_{p}$.
\end{theorem}

and this show that the problem of syzygies of the Veronese embedding is
 connected to the following problem: 
 let $Y$ be a smooth complex projective variety and $M$  a line bundle 
on $Y$;  if $M$ satisfies Property $ N_{p}$ then for which $k$
the bundle $M^{d}$ satisfies Property $N_{k}$? (see also \cite{Ru2}).

\smallskip

It seems difficult to have some result on syzygies of Veronese embedding 
under the diagonal $d=p$, especially 
a result holding for ${\bf P}^n$ for every $n$.
Here   we prove:

\begin{theorem} \label{N2}
The line bundle  $ {\cal O}_{{\bf P}^{n}}(3)$ satisfies Property $N_{4}$
 $\forall n$.
\end{theorem}

Perhaps the technique used here to prove Thm. \ref{N2} may be useful to 
solve some other open case of syzygies of Veronese embeddings.

To prove Thm. \ref{N2} we  prove also:

\begin{proposition} \label{vari}
 Let $d, n  \in {\bf N}$;  we have that  
 ${\cal O}_{{\bf P}^n}(d)$   satisfies $N_p$  $\forall n \geq  p$ 
if and only if  ${\cal O}_{{\bf P}^{p}
}(d)$  satisfies $N_p$.
\end{proposition}

\section{Proof of Proposition \ref{vari}}

{\em Proof of Prop. \ref{vari}.}
 Let   $ L $ be a very ample  line bundle on a smooth
complex projective variety $Y$. We recall from \cite{Green2} that $L$
satisfies $N_p$ iff $$ (Tor^{S(L)}_p (G(L),{\bf C}) 
)_{p+q}= 0 \;\; \forall q \geq 2$$ 
(see Introduction for the notation) and
$ (Tor^{S(L)}_p (G(L),{\bf C}) )_{p+q}$ is equal to the homology of
the Koszul complex
$$\wedge^{p+1} H^0 (L)  \otimes H^0(L^{q-1})  \rightarrow 
\wedge^{p} H^0 (L)  \otimes H^0(L^{q})
\rightarrow \wedge^{p-1} H^0 (L)  \otimes H^0(L^{q+1})$$

Now let $Y=  {\bf P}(V)$ be a projective space and
 $L={\cal O}_{{\bf P}(V)}(d)$. In our case   
$ (Tor^{S(L)}_p (G(L),{\bf C}) )_{p+q}$ is equal to the homology of
$$\wedge^{p+1} Sym^d V \otimes Sym^{(q-1)d} V 
\stackrel{\alpha^V_{p+1,q-1}}{\longrightarrow} 
\wedge^{p} Sym^d V \otimes Sym^{qd} V 
\stackrel{\alpha^V_{p,q}}{\longrightarrow}  
\wedge^{p-1} Sym^d V \otimes Sym^{(q+1)d} V $$ 
since the maps are $GL(V)$-invariant,  $ (Tor^{S(L)}_p (G(L),{\bf C}) )_{p+q}$
is a $GL(V)$-module,
 as observed in Rem. of \S 2 of \cite{Green1} and Prop. 1.8
\cite{O-P}.
The Young diagrams of the irreducible subrepresentations of
the $GL(V)$-module 
$\otimes^{p} (Sym^d V)$ have  at most $p $ rows (see for 
instance p. 79 \cite{F-H}), thus
the Young diagrams of the irreducible subrepresentations of 
$\wedge^{p} (Sym^d V)$ have  at most $p $ rows 
and then the ones  of $\wedge^{p} Sym^d V \otimes Sym^{qd} V $ 
have at most $p  +1 $ rows .Thus the Young diagrams of the
 irreducible subrepresentation of $ (Tor^{S(L)}_p (G(L),{\bf C}) )_{p+q}$ 
have at most  $p  +1 $ rows and these Young diagrams don't depend on $V$, 
in fact:

 by Littlewood-Richardson's rule we can write 
$\wedge^{p} Sym^d V \otimes Sym^{qd} V = \oplus_{\lambda \in A_{p,d,q}}
S^{\lambda} V $ where $A_{p,d,q}$ is a subset of the set of the partitions 
of $ pd+ qd$ and does not depend $V$; 
we want to show that the Young diagrams of the irreducible
 subrepresentations of
$Ker (\alpha^V_{p,q})$ and of $Im (\alpha^V_{p+1,q-1})$
 don't depend on $V$; let  $V$ and $W$ be two  vector spaces; suppose for 
instance  $dim (V) \geq dim(W) $, then there exists an injective map 
$W \rightarrow V$ and 
we have the following commutative diagram:
{\scriptsize
$$\begin{array}{ccccc} 
\wedge^{p+1} Sym^d W \otimes Sym^{(q-1)d} W & 
\stackrel{\alpha_{p+1,q-1}^W}{\longrightarrow} & 
\wedge^{p} Sym^d W \otimes Sym^{qd} W  &  
\stackrel{\alpha_{p,q}^W}{\longrightarrow} &
\wedge^{p-1} Sym^d W \otimes Sym^{(q+1)d} W
\\  \downarrow & & \downarrow & & \downarrow
\\
\wedge^{p+1} Sym^d V \otimes Sym^{(q-1)d} V & 
\stackrel{\alpha_{p+1,q-1}^V}{\longrightarrow} & 
\wedge^{p} Sym^d V \otimes Sym^{qd} V  & 
\stackrel{\alpha_{p,q}^V}{\longrightarrow} &
\wedge^{p-1} Sym^d V \otimes Sym^{(q+1)d} V 
\end{array}$$ }
that can be written as  ($S^{\lambda} $ denotes the Schur functor associated
 to $\lambda$):
$$\begin{array}{ccccc} 
\oplus_{\lambda \in A_{p+1,d,q-1}}S^{\lambda} W
 & \rightarrow & \oplus_{\lambda \in A_{p,d,q}}  S^{\lambda} W
 & \rightarrow & \oplus_{\lambda \in A_{p-1,d,q+1}} S^{\lambda} W
\\  \downarrow & & \downarrow & & \downarrow \\
 \oplus_{\lambda \in A_{p+1,d,q-1}}S^{\lambda} V
 & \rightarrow & \oplus_{\lambda \in A_{p,d,q}}  S^{\lambda} V
 & \rightarrow & \oplus_{\lambda \in A_{p-1,d,q+1}} S^{\lambda} V
\end{array}$$ 
If $\lambda \in A_{p,d,q} - A_{p-1,d,q+1}$ (which is a set not depending
 on $V$) then  $S^{\lambda} V$ is in $Ker (\alpha^V_{p,q})$;
besides if $\lambda \in A_{p,d,q} \cap A_{p-1,d,q+1}$,
 the map $S^{\lambda} W \rightarrow S^{\lambda} W$ (which
 can be only a multiple of identity by Schur Lemma) induced by  
$\alpha_{p,q}^W$ is nonzero iff
 the corresponding map  $S^{\lambda} V \rightarrow
 S^{\lambda} V$ induced by 
$\alpha^V_{p,q}$ is nonzero; thus the  Young diagrams of the irreducible
 subrepresentations of
$Ker (\alpha^V_{p,q})$  don't depend on $V $ and 
analogously for  $Im  (\alpha^V_{p+1,q-1}) $;
thus $Ker(\alpha^V_{p,q})/Im(\alpha^V_{p+1,q-1})= 
\oplus_{\lambda \in A'_{p,d,q}}  S^{\lambda} V$ for some subset 
$ A'_{p,d,q} $ of  $A_{p,d,q}$ not depending on $V$.

Since  the Young diagrams of the
 irreducible subrepresentations of $ (Tor^{S(L)}_p (G(L),{\bf C}) )_{p+q}$ 
have at most  $p  +1 $ rows and these Young diagrams don't depend on $V$, 
we have that if these representations  are zero for $dim(V) = 
p + 1 $  they are zero also for $dim(V) \geq p  +1 $.
\hfill \includegraphics[scale=0.9]{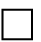}

\section{Recalls on syzygies of toric ideals}

We  recall   some facts on toric ideals from \cite{St}.
Let $k \in {\bf N}$. Let $A = \{a_1,..., a_m \} \subset {\bf Z}^k$.
The toric ideal  ${\cal I}_A$ is defined as  the  ideal 
in ${\bf C}[x_1, ..., x_m]$
generated as vector space by the
binomials 
$$x_1^{u_1}... x_m^{u_m}
- x_1^{v_1}... x_m^{v_m}$$
 for $(u_1,..., u_m),(v_1,.., v_m) \in {\bf N}^m$, 
with $\sum_{i=1,...,m} u_i a_i = \sum_{i=1,...,m} v_i a_i $. 
           
We have that ${\cal I}_A$ is homogeneous iff $\exists \;
\omega \in {\bf Q}^k$ s.t. $\omega  \cdot a_i= 1$ $\forall i = 1,..., m$;
the rings ${\bf C}[x_1, ..., x_m]$ and 
  ${\bf C}[x_1, ..., x_m]/{\cal I}_A$ 
are multigraded by ${\bf N}A$ via $\deg x_i =a_i$;  
the element $x_1^{u_1}...x_m^{u_m}$ 
has multidegree $ b =  \sum_i u_i a_i \in {\bf N}A$ and degree 
$\sum_i  u_i = b \cdot \omega$; we define $\deg b = b \cdot \omega$.

   For each $b \in {\bf N}A$, let 
$\Delta_b$ be the simplicial complex (see \cite{Sp})
on the set $A$  defined   as follows:
 $$ \Delta_b = \{ \langle F \rangle 
 | \; F \subset A : b - \sum_{a \in F} 
a  \in {\bf N}A \}$$
             

The following theorem studies the syzygies of the
ideal ${\cal I}_A$; it was proved by Campillo and Marijuan for $k=1$ 
in \cite{C-M} and by Campillo and Pison for general $k$ and $j=0$ in 
\cite{C-P}; the following more general statement is due to Sturmfels  
(Theorem 12.12 p.120 in \cite{St}).

\begin{theorem} \label{CPS}
({\rm see \cite{St}, \cite{C-M}, \cite{C-P}}).
Let $A=\{a_1,...,a_m\} \subset {\bf N}^k$
 and ${\cal I}_A$ be the associated 
toric ideal.  Let $ 0 \rightarrow E_n \rightarrow ...\rightarrow
 E_1 \rightarrow E_0 \rightarrow G \rightarrow 0$ be a 
minimal free resolution of $G= {\bf C}[x_1,..., x_m]/ {\cal I}_A$
on ${\bf C}[x_1,..., x_m]$.
Each of the generators of $E_j$ has a unique multidegree. 
The number of the generators of 
multidegree $b \in {\bf N}A$ of  $E_{j+1}$  equals the rank of the 
$j$-th
reduced homology group $\tilde{H}_j(\Delta_b, {\bf C}) $. 
\end{theorem}

\section{Proof of Theorem \ref{N2}}

\begin{notation}
$\bullet$ homologous means homologous in the reduced homology.

$\bullet$ $\sim_A$ means homologous in $A$, i.e. $\gamma \sim_A \gamma'$
means that $\exists \, \beta$ chain in $A$ s.t. $\partial \beta = 
\gamma - \gamma'$.

$\bullet$ $e_i$ denotes the $i$-th 
element of the canonical basis of ${\bf R}^n$.

$\bullet $ The symbol $\ast $ denotes the joining.

$\bullet $ For any  $v \in {\bf R}^n$ 
$v_i$ denotes the $i$-th coordinate, that is the lower index denotes the 
coordinate. 
\end{notation}

If we take $A =A_{d,n}= \{ 
\tiny{ \left( \hspace{-0.2cm}
\begin{array}{c} x_1 \\. \\. \\
x_{n+1} \end{array}  \hspace{-0.2cm} \right) }
| \sum x_i = d\;\;  x_i \in {\bf N}\}$, 
     we have that ${\cal I}_{A_{n,d}}$ is the ideal of the  
embedding of ${\bf P}^n$ by  ${\cal O}(d)$.
In this case  $\omega =\omega_d = \frac{1}{d} (1, ..., 1)$.

Let $b \in {\bf N} A_{d,n}$; we have that  $\deg b=(= b \cdot \omega) =k$ iff
$b$ is the sum of $k$ (not necessarily distinct) elements of $A_{d,n}$.
Observe that a simplex $S$ 
 with vertices in $A_{d,n}$ is a simplex of $\Delta_b$
iff  the sum $s$ of the vertices  of $S$ is s.t.  $s_i \leq b_i$ 
$\forall i =1,....,n+1$.
We generalize the definition of the simplicial complex 
$\Delta_b$ in \S 3   in the following way:

\begin{notation} Let $v \in {\bf N}^{n+1}$. Let $\Delta_v$ be  
the following simplicial complex: a
simplex $S$ with vertices in $A_{d,n}$ is a simplex of $\Delta_v$
iff  the sum $s$ of the vertices  of $S$ is s.t.  $s_i \leq v_i$ 
$\forall i =1,....,n+1$.
\end{notation}

The main points of the  proof of the Thm. \ref{N2}
  are Propositions \ref{step1}
 and \ref{step2}.

\begin{notation} Let $d,n  \in {\bf N}$ and $b \in {\bf N} A_{d,n}$. Let $X_b$ 
be the following simplicial complex on $A_{d,n}$: 
$$
X_b := \Delta_b \cup \Delta_{
\tiny{
\left( \hspace{-0.2cm} \begin{array}{c}
b_1 -1\\
b_2 +1 \\
b_3 \\
. \\
.\\
.
\end{array}
 \hspace{-0.2cm}
\right)}} \cup .... \cup \Delta_{
\tiny{
\left( \hspace{-0.2cm} \begin{array}{c}
0\\
b_2 + b_1 \\
b_3 \\
. \\
.\\
.
\end{array}
 \hspace{-0.2cm}
\right)}} $$ (in the obvious sense that a simplex with vertices in
 $A_{d,n}$  is a simplex of $X_b$ iff  it is a simplex of 
$\Delta_{b - k e_1 +k e_2}$ for some $k \in \{0,..., b_1\}$)
\end{notation}

\begin{proposition} \label{step1} Let $d,n, p \in {\bf N}$.
Let $b \in {\bf N} A_{d,n}$ with $deg(b) \geq p+2$.
Let $\gamma $ be a $(p-1)$-cycle in $\Delta_b$.
If the following conditions hold

a) $H_{p-3}(\Delta_{c- e_1  }) =0 $
 $\forall c  \in {\bf N} A_{d,n}$ with 
  $deg(c) \geq p+1$  

b) ${\cal O}(d)$ satisfies Property $N_{p-1}$,

then  $\exists \gamma'$ cycle  in $\Delta_{
\tiny{
\left( \hspace{-0.2cm} \begin{array}{c}
0\\
b_1+b_2 \\
b_3 \\
. \\
.\\
.
\end{array}
 \hspace{-0.2cm}
\right)}}$ s.t. $\gamma \sim_{X_b} \gamma'$. 
\end{proposition}

\begin{proposition} \label{step2}  
Let $b \in {\bf N} A_{3,4}$  with $deg(b) \geq 6$.
Let $\gamma $ be a $3$-cycle in $\Delta_b$. If   $\gamma \sim_{X_b} 0$
 then $ \gamma \sim_{\Delta_b} 0 $.
\end{proposition}

\begin{lemma} \label{H1(4+1/3)}  Let $d,n  \in {\bf N}$ and
 $g \in {\bf N} A_{d,n}$ with $deg (g) \geq 4$.
 Then $H_1(\Delta_{g+v})=0 $ $\forall v \in {\bf N^{n+1}}$.
\end{lemma}

\smallskip

We show now how  Thm. \ref{N2} follows  from Propositions
 \ref{step1} and  \ref{step2} and  Lemma \ref{H1(4+1/3)}.

\bigskip

{\it Proof of Thm. \ref{N2}}. By Prop. \ref{vari}, 
  it is sufficient to prove our statement when $n=4$.
By Thm. \ref{CPS}, the bundle 
 ${\cal O}_{{\bf P}^n}(d)$  satisfies $N_p$ iff 
 $H_{q-1}(\Delta_b)=0$  $\forall b \in {\bf N}A_{d,n}$ with 
$\deg b  \geq q+2$ $\forall q \leq p$;
in particular,  in order to prove that  
 ${\cal O}_{{\bf P}^n}(3)$ 
 satisfies $N_4$, we have to prove 
that $H_{3}(\Delta_b)=0$  $ \forall b \in {\bf N}A_{3,n}$ with 
$\deg b  \geq 6$. Thus
let $b \in {\bf N} A_{3,4}$ with $deg(b) \geq 6$.
Let $\gamma $ be a $3$-cycle in $\Delta_b$. We want to prove  
$\gamma \sim_{\Delta_b} 0$.  
 By   Prop. \ref{step1}, $ \exists \gamma'$ 
cycle   in $\Delta_{
\tiny{
\left( \hspace{-0.2cm} \begin{array}{c}
0\\
b_1+b_2 \\
b_3 \\
. \\
.\\
.
\end{array}
 \hspace{-0.2cm}
\right)}}$ s.t.  $\gamma \sim_{X_b} \gamma '$ 
(the assumption of Prop. \ref{step1} 
in our case  is true by Lemma \ref{H1(4+1/3)}).  
We have $H_3 (\Delta_{
\tiny{
\left( \hspace{-0.2cm} \begin{array}{c}
0\\
b_1+b_2 \\
b_3 \\
. \\
.\\
.
\end{array}
 \hspace{-0.2cm}
\right)}})= H_3 (\Delta_{
\tiny{
\left( \hspace{-0.2cm} \begin{array}{c}
b_1+b_2 \\
b_3 \\
. \\
.\\
.
\end{array}
 \hspace{-0.2cm}
\right)}})
 =0$, where the last equality holds since 
 ${\cal O}_{{\bf P}^3}(3)$  satisfies $N_4$ (by Prop. \ref{OP3},
 but it can be proved also directly).
Thus   $\gamma' \sim 0 $ in $\Delta_{
\tiny{
\left( \hspace{-0.2cm} \begin{array}{c}
0\\
b_1+b_2 \\
b_3 \\
. \\
.\\
.
\end{array}
 \hspace{-0.2cm}
\right)}}$.  
Thus  $\gamma \sim_{X_b} 0 $ and then
 $\gamma \sim_{\Delta_b} 0$
by Prop. \ref{step2}.
\hfill \includegraphics[scale=0.9]{box}

\bigskip

\bigskip

Now we will prove  Propositions
 \ref{step1} and  \ref{step2} and  Lemma \ref{H1(4+1/3)}.

\begin{notation}
Let $b \in {\bf N} A_{d,n}$
Let $\gamma$ be a $(p-1)$-cycle in $X_b$.

 For every vertex $a$  in $\gamma$, let
${\cal S}_{a,\gamma} $ be  the set of simplexes of $\gamma$ with vertex  $a$ 
and $\mu_{a, \gamma}$ be the $(p-2)$-cycle s.t. $a \ast \mu_{a, \gamma} 
= \sum_{\tau \in {\cal S}_{a,\gamma}} \tau$.  
For $\tilde{a} \in A_{d,n}$, let 
$$\alpha_{a, {\tilde a}, \gamma} = (a- \tilde{a}) \ast  \mu_{a, \gamma}$$
\end{notation}

\bigskip
 
{\it Proof of Prop. \ref{step1}}.
We order in some way the (finite) vertices of  
$\gamma$ with first coordinate $ \neq 0$: $a^1,..., a^r$.
Let $\tilde{a}^{j}= \tiny{
\left( \hspace{-0.2cm} \begin{array}{c}
0 \\
a^{j}_2 +a^{j}_1\\
a^{j}_3 \\
a^{j}_4 \\
.\\
.
\end{array}
 \hspace{-0.2cm}
\right)}$ for $j=1,..., r$.

Obviously  $\alpha_{a^1,\tilde{a}^1 , \gamma} \sim_{X_b} 0 $,
because  $\mu_{a^1, \gamma} $
 is in $\Delta_{b-a^1}$ and 
$H_{p-1}((a^1 - \tilde{a}^1) \ast  \Delta_{b-a^1}) =
H_{p-2}(\Delta_{b-a^1}) =0$ (since    
${\cal O}(d)$ satisfies Property $N_{p-1}$).
Thus $\gamma_1 := \gamma + \alpha_{a^1, \tilde{a^1},
 \gamma}$ is homologous to $\gamma$ in $X_b$.

We define by induction $\gamma_j := \gamma_{j-1} + \alpha_{a^j,
\tilde{a}^j,  \gamma_{j-1}}$ for $j=2,..., r$.
We want to  prove $\gamma_r \sim_{X_b} 0$; to  prove this, we
 prove  $\alpha_{a^j, \tilde{a^j}, \gamma_{j-1}} \sim_{X_b} 0$
for $j=2,..., r$.

Observe that 
 $\mu_{a^j, \gamma_{j-1} } $ is in
$$\Delta_{b-a^j} \cup \Delta_{
\tiny{
\left( \hspace{-0.2cm} \begin{array}{c}
(b-a^j)_1 -1\\
(b-a^j)_2 +1 \\
(b-a^j)_3 \\
. \\
.\\
.
\end{array}
 \hspace{-0.2cm}
\right)}} \cup .... \cup \Delta_{
\tiny{
\left( \hspace{-0.2cm} \begin{array}{c}
0\\
(b-a^j)_1+(b-a^j)_2 \\
(b-a^j)_3 \\
. \\
.\\
.
\end{array}
 \hspace{-0.2cm}
\right)}}$$

We can find some cycles $\theta_{\varepsilon}$ 
in  $\Delta_{b-a^j - \varepsilon e_1 + \varepsilon e_2}$ for 
 $\varepsilon \in \{0, ..., (b-a^j)_1\}$ s.t. 
 $\mu_{a^j, \gamma_{j-1} }= \sum_{\varepsilon \in \{0, ..., (b-a^j)_1\}}
 \theta_{\varepsilon}$, in fact: 
let $\sigma_0$ be the sum of the simplexes of $\mu_{a^j, \gamma_{j-1} }$ in 
$\Delta_{b-a^j}$ and not in   $\Delta_{b-a^j - e_1 }$; 
$\partial \sigma_0 $ is in  $\Delta_{b-a^j -  e_1 }$ and 
since $H_{p-3}(\Delta_{c- e_1}) =0 $  $\forall c $ with $ deg(c) \geq p+1$
then
$\exists \sigma_0' $ in $\Delta_{b-a^j -  e_1 }$ s.t. 
$\partial \sigma_0' = \partial \sigma_0 $; 
let $ \theta_0= \sigma_0 - \sigma_0'$; now 
$\mu_{a^j, \gamma_{j-1} } - \theta_0$ is in $$ \Delta_{\tiny{
\left( \hspace{-0.2cm} \begin{array}{c}
(b-a^j)_1 -1\\
(b-a^j)_2 +1 \\
(b-a^j)_3 \\
. \\
.\\
.
\end{array}
 \hspace{-0.2cm}
\right)}} \cup .... \cup \Delta_{
\tiny{
\left( \hspace{-0.2cm} \begin{array}{c}
0\\
(b-a^j)_1+(b-a^j)_2 \\
(b-a^j)_3 \\
. \\
.\\
.
\end{array}
 \hspace{-0.2cm}
\right)}}$$ 
and we can go on analogously: let  $\sigma_1$ 
be the sum of the simplexes of $\mu_{a^j, \gamma_{j-1} }- \theta_0$ in 
$\Delta_{b-a^j-e_1 +e_2}$ and not in   $\Delta_{b-a^j - 2 e_1 }$....

Since $deg(b-a^j -\varepsilon e_1 + \varepsilon e_2 ) \geq p+1$
 and   ${\cal O}(d)$ satisfies $N_{p-1}$ we have
$H_{p-2} (\Delta_{b-a^j -\varepsilon e_1 + \varepsilon e_2} )=0$   thus
$ \theta_{\varepsilon} \sim 0$ in  $\Delta_{
b-a^j -\varepsilon e_1 + \varepsilon e_2}$ 
and then $(a^j -\tilde{a^j}) \ast \theta_{\varepsilon} \sim_{X_b} 0$;
therefore $\alpha_{a^j, \tilde{a}^j, \gamma_{j-1}} \sim_{X_b} 0$.

Thus we can define $\gamma' = \gamma_r$: 
 $\gamma' $ is  in  $\Delta_{
\tiny{
\left( \hspace{-0.2cm} \begin{array}{c}
0\\
b_1+b_2 \\
b_3 \\
. \\
.\\
.
\end{array}
 \hspace{-0.2cm}
\right)}}$ and $ \gamma' \sim_{X_b} \gamma$.
\hfill \includegraphics[scale=0.9]{box}

\bigskip

{\it Proof of Lemma \ref{H1(4+1/3)}.} 
By induction on the sum of the coefficients of $v$. If $v=0$ the statement is 
true  since  ${\cal O}_{{\bf P}^n} (d)$
satisfies $N_2$ $\forall n$ $\forall d \geq 2$. 

Let $h =g+v$.
Suppose   $H_1(\Delta_h) = 0 $; we want to show $H_1 (\Delta_{h+e_j}) =0$;
 it is sufficient to prove that every $1$-cycle  $\gamma $ in $\Delta_{h+e_j}$ 
is homologous in  $ \Delta_{h +e_j}$
to some 1-cycle   in  $ \Delta_{h}$.
Let $ \langle x, a \rangle  $ be a simplex of $\gamma $  not in  $ \Delta_{h}$.
Thus or $x_j > 0$ either $a_j >0$, say for instance 
$a_j >0$. Let $y$ be a vertex in $\mu_{a,\gamma}$ with $y \neq x$. 
Let $i$ be s.t. $x_i +y_i +a_i < h_i$
(such an $i$ exists because $deg(g) \geq 4$). Let $\tilde{a} = a +e_i -e_j$ and
 $ \alpha =\langle x,a  \rangle  + \langle a, y \rangle  
+\langle y, \tilde{a} \rangle  + \langle \tilde{a}, x  \rangle $; then $\alpha 
\sim_{\Delta_{h+e_j}} 0 $ because $\alpha$ 
 is in $ (a -\tilde{a}) \ast \Delta_{h + e_j  -a - e_i}  $  and 
$H_1((a - \tilde{a}) \ast \Delta_{h + e_j  -a - e_i}  ) = \tilde{H}_0 (
\Delta_{h + e_j  -a - e_i}) =0$  ($\Delta_{h + e_j  -a - e_i}$ is connected 
since if $v,w $ are two vertices in $ \Delta_{h + e_j  -a - e_i}$ we can find 
$u \in A_{d,n}$
 s.t. $ \langle u,v \rangle $ and $ \langle u, w \rangle $ are in 
$\Delta_{h + e_j  -a - e_i}$ since
 $ deg (g + e_j  -a - e_i) \geq 3$).
Let $ \gamma' = \gamma + \alpha$. We have  $\gamma' \sim_{
\Delta_{h +e_j}} \gamma$ 
and the number of the simplexes of $\gamma'$ not in $\Delta_h$ 
is less than the number of simplexes of $\gamma$ not in $\Delta_h$. 
Thus, by induction on the number of simplexes of $\gamma$ not in $\Delta_h$, 
we get a $1$-cycle in $\Delta_h$ homologous  in 
$\Delta_{h+e_j}$ to $\gamma$.
\hfill \includegraphics[scale=0.9]{box}

\bigskip
Now we will prove Prop. \ref{step2}.

\begin{definition} \label{UFO} Let $d,n,k \in {\bf N}$ and
 $\beta \in {\bf N} A_{d,n}$.
We say that a  $(k-1)$-chain   $\eta$  in $\Delta_{\beta}$ 
  is a $UFO$ with axis  $\langle a^1,..., a^t \rangle $ for the
 coordinate $i$ (for short we will write   
$UFO_{t,k}^{i}(a^1,..., a^t, \Delta_{\beta})$)
if  
$$\eta = \langle a^1,..., a^t \rangle  \ast C_{\eta}$$ for some 
$ (k-t -1)$-cycle $C_{\eta}$ and
 $a^1,..., a^t $ are distinct vertices in $  \Delta_{\beta}$
 with $$(a^1 + ...+ a^t)_i = \beta_i  \;\;\;\;\;\;\;
(a^j)_i > 0 \; \; \forall j =1,...,l$$
We will denote the axis 
$\langle a^1 ,..., a^t \rangle$ by $\chi_{\eta}$.  Observe 
$\partial \eta \subset \Delta_{\beta -e_i}$.
 
(Sometimes we will omit some index when it will be obvious.)
\end{definition}

\smallskip

\hspace{5.5cm}
\includegraphics[scale=0.3]{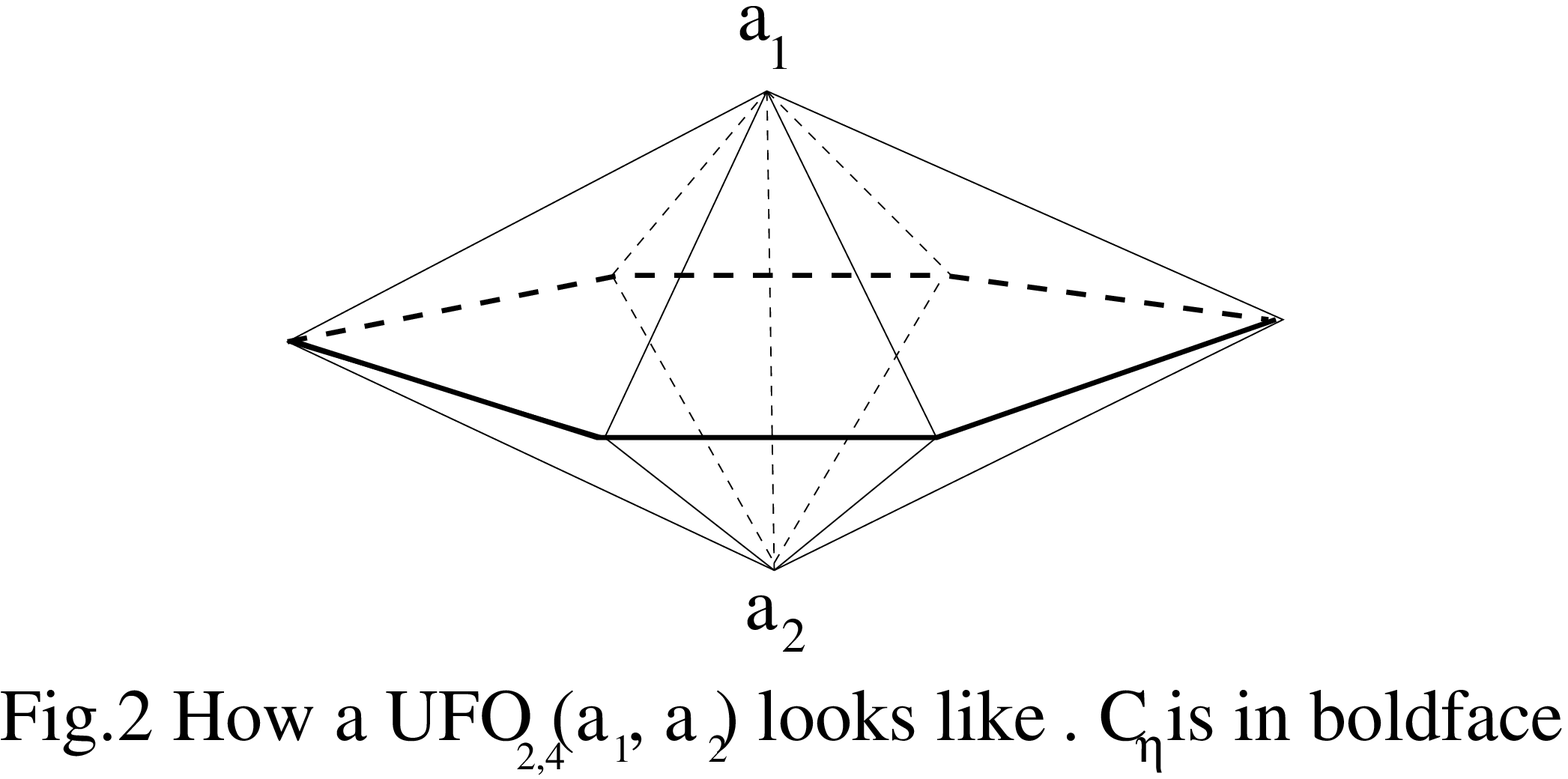}


\begin{lemma} \label{ufisemplici} {\boldmath
 ``$UFO_{p+1, p+1}$, `$UFO_{p, p+1}$, 
`$UFO_{1, p+1}$ ''}
 Let $d,n, p \in {\bf N}$ and  $r,l  \in \{1,...,n+1\}$, $r \neq l$.
 Let $\beta \in {\bf N} A_{d,n}$ with  $\deg \beta \geq p+2$.
Let $\eta $ be  a $UFO_{t, p+1}^r (a^1, ..., a^{t},\Delta_{\beta})$. 

 If $t \in \{p+1, p , 1\}$ 
then $\exists \tilde{\eta} $ $p$-chain  in $\Delta_{
\beta +e_l -e_r}$ with $\partial \tilde{\eta}= \partial \eta$.

\end{lemma}

{\it Proof.} \underline{Case $t=p+1$}
 Since $\deg \beta \geq p+2 $ and $\eta $ is a simplex with 
$p+1$ vertices (the simplex $\langle a^1,..., a^{p+1} \rangle$), then     
 $\exists x \in A_{d,n}$  s.t. $ x \ast \eta  \subset \Delta_{\beta}$.
 Since $(a^1+ ...+a^{p+1})_r = \beta_r$ then  $x_r=0$. Take  $\tilde{\eta} :=
x \ast \partial \eta $.

\underline{Case $t=p$}
 If $\deg \beta \geq p+3$ then $ \deg (\beta -a^1 ...- a^p ) \geq 3$, 
therefore $ \tilde{H}_0 (\Delta_{\beta -a^1 ...- a^p}) =0$, thus  
$\exists \gamma $ s.t. $ \partial \gamma = C_{\eta}$
 and we can take $ \tilde{\eta}=
\gamma \ast \partial  \chi_{\eta}$.

Thus we can suppose $\deg \beta = p+2$. 
Since $C_{\eta}$ is a $0$-cycle it is sufficient to prove the statement when
$C_{\eta}=P-Q$ for some $P,Q \in A_{d,n}$ with
$P=Q +e_i -e_j$  for some $i$ and $j$ (in fact  we can write $C_{\eta}$ as 
$ \sum_s (-1)^s P_s$ with $P_{s+1} $ obtained from $P_s$  by adding 
$1$ to a  coordinate   and subtracting $1$ to another coordinate). 
Let $x = \beta -a^1...-a^p -P$ and  $y =\beta -a^1...-a^p -Q$.
Since $(a^1 + ...+ a^p)_r = \beta_r$ we  have  $x_r =y_r =0$.

Suppose first  $i, j \not \in \{l,r\} $.
Let $z = x +e_l -e_j$.
The chain  $\tilde{\eta } = z  \ast \partial \eta$ is in $ \Delta_{\beta +e_l
-e_r}$ and $\partial \tilde{\eta} = \partial \eta$.

Suppose now   $j \not \in \{l ,r\}$  and $i=l$.
Then  $\tilde{\eta } := x  \ast \partial \eta$ is in $ \Delta_{
\beta +e_l -e_r}$ and $\partial \tilde{\eta} = \partial \eta$.

\underline{Case $t=1$} Let $ \tilde{a^1} = a^1 +e_l -e_r$;
take  $\tilde{\eta} = \tilde{a^1}  \ast C_{\eta }$.
 \hfill \includegraphics[scale=0.9]{box}

\begin{lemma} \label{subc}  Let $d,n, p \in {\bf N}$ and
$r \in \{1,...,n+1\}$.
Let $\beta \in {\bf N} A_{d,n}$ and $\deg \beta \geq  p+2$.  
Let $\eta $ be  a $UFO_{t, p+1}^r (a^1,..., a^t, \Delta_{\beta })$. 
If $C_{\eta} = \partial \sigma$ where $\sigma $ is a simplex 
(with $p+2 -t $   vertices) in
 $\Delta_{\beta - a^1 -... -a^t} $, then $\exists \tilde{\eta}$ $p$-chain 
in $\Delta_{\beta  -e_r}$ with $\partial \tilde{\eta}= \partial \eta$.
\end{lemma}

{\it Proof.} Since $\eta = \chi_{\eta}   \ast
C_{\eta} =  \chi_{\eta}  \ast   \partial \sigma $, then
 $\partial \eta =  \partial \chi_{\eta}  * \partial \sigma = \partial (
\partial \chi_{\eta}  *  \sigma)$. Take $\tilde{\eta} 
= \partial \chi_{\eta}  *  \sigma$. 
\hfill \includegraphics[scale=0.9]{box}

\smallskip

\begin{lemma} \label{ufo3} {\boldmath ``$UFO_{3,5}$''} 
Let $\beta \in {\bf N} A_{3,4}$ and $\deg \beta \geq 6$.
Let $r,l \in \{1,...,5\}$, $r \neq l$.  
Let $\eta $ be  a $UFO_{3, 5}^r (a^1, a^2, a^3,\Delta_{\beta})$. 
 Then   $\exists \tilde{\eta}$ $4$-chain in $\Delta_{
\beta +e_l -e_r}$ with $\partial \tilde{\eta}= \partial \eta$.
\end{lemma}



To prove this lemma, some  sublemmas are necessary.

\begin{sublemma} \label{suba} The statement of Lemma \ref{ufo3} is true if  
 $\eta$ is in $\Delta_{\beta -e_l}$ and
 $(a^1 + a^2 +a^3 )_l \leq \beta_l -2$.
\end{sublemma}

{\it Proof.}  
Let $\tilde{a^i} = a^i -e_r +e_l$ for $i=1,2,3$.
Observe that $- \langle  \tilde{ a^1}, a^2, a^3 \rangle  \ast C_{\eta}$
 is in $\Delta_{\beta - e_r}$.
Thus we can sum it to $\eta $ and  prove that the border of the sum is 
homologous to $0$ in $\Delta_{\beta -e_r + e_l}$. The sum is 
$(\langle a^1, a^2, a^3 \rangle - \langle \tilde{a^1} , a^2, a^3 \rangle )
 \ast C_{\eta}$ and has the same border of  the sum of the six chains:
$$ \begin{array}{lcccr}
\sigma_1= \langle a^1, a^2, \tilde{a^3}  \rangle  \ast C_{\eta },& &
 \sigma_2= \langle a^1,  \tilde{a^3} , 
\tilde{a^2} \rangle  \ast C_{\eta }, & &
 \sigma_3 = \langle a^1, \tilde{a^2} , a^3 \rangle  \ast C_{\eta } \\
 \sigma_4= \langle \tilde{a^1} ,  \tilde{a^3}, a^2  \rangle  \ast C_{\eta },
& & \sigma_5= \langle \tilde{a^1 } ,  \tilde{a^2}  , \tilde{a^3}  
\rangle  \ast C_{\eta },
& & \sigma_6= \langle \tilde{a^1} , a^3, \tilde{a^2}  \rangle  
\ast C_{\eta } 
\end{array}$$
Observe that $\sigma_1$, $\sigma_2$, $\sigma_3$,
$ \sigma_4$ and $\sigma_6$ are in $\Delta_{\beta - e_r + e_l}$
since $\eta$ is in $\Delta_{\beta -e_l}$; besides 
$\sigma_5 $ is in $\Delta_{\beta -3e_r + 2e_l }$ and its border is
 in  $\Delta_{\beta -3e_r +e_l}$;
$\sigma_5$  is  not a
$UFO^l_{3, 5}(\Delta_{\beta -2e_r +2e_l})$ since 
   $(\tilde{a}^1 + \tilde{a}^2 + \tilde{a}^3 )_l  
=(a^1 + a^2 +a^3 )_l+3 \leq \beta_l +1$; besides  
$\tilde{a^1 }_l ,  \tilde{a^2}_l  , \tilde{a^3}_l  >  0$;
 thus $\sigma_5$ is a 
sum of $UFO^l_{t,5}(\Delta_{\beta -2e_r +2e_l})$  with $ t = 4,5$; then   
  $\partial \sigma_5 \sim_{\Delta_{\beta -e_r +e_l}} 0$ by Lemma 
\ref{ufisemplici}
 (applied with $r$ equal to $l$  and $l$ equal to $r$).
 \hfill \includegraphics[scale=0.9]{box}
{\it Sl. \ref{suba}}

\begin{sublemma} \label{subb}  
The statement of Lemma \ref{ufo3} is true if  $\eta $ is in 
$\Delta_{\beta - e_j -e_i} $  for some $i$ and $j$ 
  ($i$ and $j$ possibly equal and  possibly equal 
to $l$).
\end{sublemma}

{\it Proof.} Let 
$\tilde{a^1} = a^1 -e_r +e_l, \;\;\;\;\;
\tilde{a^2}= a^2 -e_r +e_i,\;\;
\;\;\;\tilde{a^3}= a^3 -e_r +e_j$.

Observe that
 $- \langle  \tilde{ a^1}, a^2, a^3 \rangle  \ast C_{\eta}$ is
 in $\Delta_{\beta - e_r}$.
Thus we can sum it to $\eta $ and prove that the border of the sum is 
homologous to $0$ in $\Delta_{\beta -e_r + e_l}$. The sum is 
$(\langle a^1, a^2, a^3 \rangle   - \langle \tilde{a^1} , a^2, a^3 \rangle )
 \ast C_{\eta}$ and has the same border of   the sum of the six chains:
$$ \begin{array}{lcccr}
\sigma_1= \langle a^1, a^2, \tilde{a^3}  \rangle  \ast C_{\eta}, & & 
 \sigma_2= \langle a^1,  \tilde{a^3} , \tilde{a^2} \rangle  \ast C_{\eta},
& & 
\sigma_3 = \langle a^1, \tilde{a^2} , a^3 \rangle  \ast C_{\eta} \\
 \sigma_4= \langle \tilde{a^1} ,  \tilde{a^3}, a^2  \rangle  \ast C_{\eta},
& & 
 \sigma_5= \langle \tilde{a^1 } ,  \tilde{a^2}  , \tilde{a^3}  \rangle  
\ast C_{\eta}, & & 
 \sigma_6= \langle \tilde{a^1} , a^3, \tilde{a^2}  \rangle  \ast C_{\eta}
\end{array} $$
Since  $\sigma_l$ for $l = 1,...,6$ 
are in $\Delta_{\beta - e_r + e_l}$,  we conclude.

(Observe that the same proof works also if $i$ or $j $ are equal to
 $l$ and if $i=j$).
\hfill \includegraphics[scale=0.9]{box}
 {\it  Sl. \ref{subb}}

\begin{sublemma} \label{subd} 
In the hyptheses  of Lemma \ref{ufo3} and  if  $\exists i $ s.t.
$(\beta -a^1 -a^2 -a^3)_i \geq 2$, 
then $\exists \tilde{\eta}$  $UFO_{3, 5}^r (a^1, a^2, a^3, \Delta_{\beta})$ 
contained in $\Delta_{\beta -e_i}$ s.t. 
$\partial \eta  \sim_{\Delta_{\beta -e_r +e_l}} \partial \tilde{\eta}$
($i$ and $l$ possibly equal) (observe the axis of $\eta$ is equal to axis of 
$\tilde{\eta}$). 
\end{sublemma}

{\it Proof.}  Observe that  $C_{\eta}$ is a  1-cycle.

$\bullet$ First  we prove that $\partial \eta \sim_{
 \Delta_{\beta -e_r + e_l}}   \partial
 \eta' $ for some $\eta'$  $UFO_{3, 5}^r (a^1, a^2, a^3, \Delta_{\beta})$
s.t. $\not \exists V$  vertex  in $C_{\eta' }$  
with $(V)_i =(\beta -a^1 -a^2 -a^3)_i$.
 
Suppose $V$ is a vertex of $C_{\eta}$ with $(V)_i =(\beta -a^1 -a^2 -a^3)_i$. 
Let $A$ and $B$ be two distinct  
vertices of $C_{\eta}$  s.t.
 $\langle A, V \rangle  + \langle V,B \rangle $ is in $C_{\eta}$.   
We define $A \cap B$  the element of ${\bf N}^{5}$ 
 s.t.  $(A\cap B)_j = min \{A_j, B_j\}$ $\forall j= 1,..., 5$.
By Lemma \ref{subc} we can suppose $ A +B +V \neq \beta -a_1 -a_2 -a_3$.
If $deg(\beta) =6$ up to changing $A$ with 
$\beta - a^1 -a^2 -a^3 -A -V$ (using Lemma \ref{subc}) 
we can suppose that the sum of the coordinates of $A\cap B$ is $2$;
thus $\exists k, j $ s.t. $(V+ A)_l ,  (V+B)_l \leq (\beta - 
a^1 -a^2 -a^3 -e_k -e_j)_l$  $\forall l$ and this  is obviously true also 
if $deg(\beta) \geq 7$;
 consider the cycle $K$ having as ordered set of vertices 
$$V, A , V- e_i +B -(A \cap B) ,  e_i + A \cap B , V -e_i +A - (A \cap B), 
B ,V $$
 Every edge of $K$  is in $\Delta_{\beta - a^1 -a^2 -a^3 -e_k -e_j}$; then 
 by Sublemma \ref{subb},  
$\partial (\langle a^1, a^2, a^3 \rangle  \ast K) \sim 0$ 
in $\Delta_{\beta -e_r +e_l}$.
Since $A_i =B_i =0$
 the vertices of $K$ different from $V$ 
 have the $i$-th cooordinate $ < (\beta -a^1 
-a^2 -a^3)_i$. 

Let $\eta' = \eta + \langle a^1, a^2, a^3 \rangle  \ast K$ 
(roughly speaking we are ``replacing'' in $C_{\eta}$
$\langle A,V \rangle  +\langle V,B \rangle $  with  an opportune  chain).

$\bullet $ Now we prove  $\partial \eta' \sim_{
 \Delta_{\beta -e_r}}  \partial
 \tilde{\eta} $ for some $\tilde{\eta}$ 
 $UFO_{3, 5}^r (a^1, a^2, a^3, \Delta_{\beta})$
s.t. $\not \exists \langle F, G \rangle$ simplex 
 in $C_{\tilde{\eta} }$ 
 with  $(F+G)_i = (\beta -a^1 -a^2 -a^3)_i $, $F_i > 0$, $G_i >0 $ and
s.t. $\not \exists V$ vertex in   $C_{\tilde{\eta} }$  
with $(V)_i =(\beta -a^1 -a^2 -a^3)_i$ (thus $\tilde{\eta}$ is 
in $ \Delta_{\beta -e_i}$).

Suppose there is a simplex  $\langle F, G \rangle $ of $C_{\eta'}$
 s.t.  $(F+G)_i = (\beta -a^1 -a^2 -a^3)_i $, $F_i > 0$ and $G_i >0 $.
Let $P$ be s.t. $\langle P,F,G \rangle \subset \Delta_{\beta-a^1 -a^2 -a^3}$. 
 By Lemma \ref{subc}, $\partial (\langle a^1, a^2,a^3 \rangle  \ast
 (\langle P,F \rangle  +\langle F,G \rangle +\langle G,P \rangle ))
 \sim_{\Delta_{\beta -e_r}} 0$.
Thus  $\partial \eta' \sim_{\Delta_{\beta -e_r}} 
\partial (\langle a^1, a^2,a^3 \rangle  \ast
 (C_{\eta'}-\langle P,F \rangle  -\langle F,G \rangle - \langle G,P 
\rangle ))$ (roughly speaking we are ``replacing''  
in $C_{\eta'}$  the simplex $ \langle  F, G \rangle $ with 
 $ \langle  F , P \rangle  +    \langle  P, G \rangle $). 
Observe that $\langle  F , P \rangle  +    \langle  P, G \rangle $ 
is in $\Delta_{\beta -a^1 -a^2 -a^3 -e_i}$. 
Repeating this for every edge $\langle F, G \rangle $ of  $C_{\eta'}$
 s.t.  $(F+G)_i = (\beta -a^1 -a^2 -a^3)_i $, $F_i > 0$ and $G_i >0 $, 
we get $\tilde{\eta}$.
\hfill \includegraphics[scale=0.9]{box}
 {\it Sl. \ref{subd}}

\bigskip


{\it Proof of Lemma \ref{ufo3}.}
Observe that  $\exists i$ s.t. 
$(\beta -a^1 -a^2 -a^3)_i \geq 3$ (in fact $\sum_j  (\beta -a^1 -a^2 -a^3)_j 
=9  $ and $(\beta -a^1 -a^2 -a^3)_2=0$).

By the Sublemma \ref{subd} we can suppose that $\eta$ is in
 $\Delta_{\beta -e_i}$.
If $i=l$ we conclude by Sublemma \ref{suba}; otherwise we can see $\eta$ as
$UFO^r_{3,5}(a^1,a^2,a^3, \Delta_{\alpha})$ where $\alpha=
\beta +e_l -e _i $; since $(\alpha -a^1 -a^2 -a^3)_i = 
(\beta -a^1 -a^2 -a^3)_i -1 \geq 2$,   by Sublemma \ref{subd} we can suppose 
$\eta$ is in $\Delta_{\alpha -e_i}$ up to homology in 
$\Delta_{\alpha -e_r +e_i} = \Delta_{\beta -e_r +e_l}$ (take $l$ of Sublemma
 \ref{subd} equal to $i$); 
  by Sublemma \ref{suba} (with $l$ of Sublemma \ref{suba} equal to $i$) 
we can conclude that $\eta \sim 0 $ in 
$\Delta_{\alpha -e_r +e_i}= \Delta_{\beta +e_l -e_r}$.  
\hfill \includegraphics[scale=0.9]{box}

\smallskip

\begin{lemma} \label{ufo2}  {\boldmath ``$UFO_{2,5}$''}  
Let $\beta \in {\bf N}  A_{3,4}$ with  $\deg \beta \geq 6$.
Let $r,l \in \{1,...,5\} $, $r \neq l$.
Let $\eta $ be a $UFO_{2, 5}^r (a^1, a^2, \Delta_{\beta})$.
 Then  
 $\exists \tilde{\eta}$ $4$-chain in $\Delta_{\beta
+e_l -e_r}$ with $\partial \tilde{\eta}= \partial \eta$.
\end{lemma}

{\it Proof.} 
$\bullet $
First suppose that $\exists i$ s.t. $(a^1 +a^2)_i <  \beta_i$ 
and at least one of $(a^1)_i$, $(a^2)_i$ is  $ > 0$, say    $(a^2)_i > 0$.

Let  $\lambda := \langle \tilde{a^1}, a^2 \rangle  \ast C_{\eta}$, 
where $\tilde{a^1}= a^1 - e_r  + e_i $.
It is not a $UFO^{i}_{2,5} (\tilde{a}^1, a^2, \Delta_{\beta -e_r +e_i})$
since $(\tilde{a^1} + a^2)_i = (a^1 +a^2)_i+1 \leq \beta_i$;
 besides $(\tilde{a}^1)_i > 0$ and $(a^2)_i > 0$; thus $\lambda $ is a sum 
of $UFO_{t, 5}^i(\Delta_{\beta}-e_r +e_i)$
for $t \geq  3$ or chains in $\Delta_{b-e_r}$ and then by Lemmas 
\ref{ufisemplici}, \ref{ufo3} (with $r$ of these lemmas  equal to $i$)
 $\partial \lambda \sim_{\Delta_{\beta +e_l -e_r}} 0 $;
 thus it is sufficient to prove $\partial (\eta + \lambda) \sim_{
\Delta_{\beta +e_l -e_r}} 0$. Let $\tilde{a^2} = a^2 -e_r +e_l$.

We have that $ \partial (\eta + \lambda)= (a^1  - \tilde{a^1}) 
\ast C_{\eta}= (a^1  - \tilde{a^2}
+  \tilde{a^2} - \tilde{a^1}) \ast C_{\eta}$; we have
$ (a^1  - \tilde{a^2}) \ast C_{\eta} \sim_{\Delta_{\beta -e_r +e_l}
} 0$ 
since   $\langle a^1 , \tilde{a^2} \rangle \ast C_{\eta}$
is  in $\Delta_{\beta -e_r +e_l}$; as in the previous paragraph we can prove 
 $  (\tilde{a^2 } - \tilde{a^1}) \ast C_{\eta} 
\sim_{\Delta_{\beta- e_r +e_l}} 0 $
(observe $\langle \tilde{a^2 }, \tilde{a^1} \rangle 
 \ast C_{\eta}$ is  in $\Delta_{\beta- 2e_r +e_l +e_i}$).

$\bullet$
In the remaing cases we have that $\forall i$ 
$(a^1 +a ^2)_i $ is equal to 
$0$ or to $\beta_i$.

Let $i$ be s.t. $ (\beta -a^1 -a^2)_i \geq 4$ (such an $i$ exists because 
the indices $j$ s.t. $(a^1 + a^2 )_j \neq 0 $ are at least two, thus the 
indeces $j$ s.t.  $ (\beta -a^1 -a^2)_j \neq 0 $ are at most three).

Let $\tilde{a^1} = a^1 -e_r +e_l$ and $\tilde{a^2} = a^2 -e_r +e_i$.
 
Observe that $\langle \tilde{ a^1}, a^2 \rangle  \ast C_{\eta}$ is in
 $\Delta_{\beta - e_r +e_l}$.
Thus we can sum it to $\eta $ and to prove that the border of the sum is 
homologous to $0$ in $\Delta_{\beta -e_r + e_l}$. The sum is 
$(\langle a^1, a^2 \rangle   + \langle a^2, \tilde{a^1} \rangle ) \ast
 C_{\eta}$ and has the same 
border of $\rho := 
(\langle a^1, \tilde{a^2} \rangle
   + \langle \tilde{a^2}, \tilde{a^1} \rangle ) \ast C_{\eta}$, thus it is
 sufficient to prove $\partial 
\rho  
\sim_{\Delta_{\beta -e_r +e_l}} 0$.

Observe $\rho $ is in $ \Delta_{\beta -e_2 +e_l +e_i}$.
Besides 
$V_i < (\beta -a^1 -a^2)_i$ $\forall V$ vertex of $C_{\eta}$ (since $V$ is
in $A_{3,4}$ and $ (\beta -a^1 -a^2)_i \geq 4$) and 
 $\tilde{a^2}$ is a vertex of every simplex
of $\rho$ and $(\tilde{a^2})_i > 0$. Thus every simplex in $\rho$ 
or is in $ \Delta_{\beta -e_2 +e_l}$ either has $3$ or $4$ vertices with 
$i$-th coordinate $>0$.
Observe that for 
any $\langle V^1, V^2, V^3  \rangle $  simplex of $C_{\eta}$ s.t.
 $(V^h)_i > 0 $ $ \forall h=1,2,3$
 and $(V^1 + V^2 +V^3 )_i = (\beta -a^1 -a^2)_i$, 
$\exists j$ s.t.   $\langle V^1 , V^2  ,V^3 \rangle \subset
  \Delta_{\beta -a^1 -a^2 -e_j} $ and 
also if  $\langle V^1, V^2, x  \rangle $  and 
  $\langle V^1, V^2, y  \rangle $ are simplexes 
of $C_{\eta}$ s.t.
 $(V^h)_i > 0 $ $ \forall h=1,2$ and $(V^1 + V^2  )_i = (\beta -a^1 -a^2)_i$, 
$\exists j$ s.t. $\langle V^1, V^2, x  \rangle $  
and   $\langle V^1, V^2, y  \rangle $ are in
 $ \Delta_{\beta -a^1 -a^2 -e_j }$  (after supposing by Lemma \ref{subc} that
$V^1 +V^2 +x+y \neq \beta-a^1 -a^2$).

 Thus   $\rho$  is a sum of  chains in $\Delta_{\beta - e_r }$ and 
$UFO^{i}_{t,5}(\Delta_{\beta - e_r +e_l + e_i -e_j})$
for $t=3,4$ and for some $j$; thus by Lemmas \ref{ufo3} and \ref{ufisemplici} 
their borders  are homologous to $0$ in $\Delta_{\beta -e_r +e_l}$.
\hfill \includegraphics[scale=0.9]{box}

\begin{corollary} \label{fine} Let $\beta \in {\bf N} A_{3,4}$ and 
$\deg \beta \geq 6$.
If $\eta$ 
is a $4$-chain  in $\Delta_{\beta}$ with 
 $\partial \eta $ in $\Delta_{
\beta  -e_2}$,
then  $\exists \tilde{\eta}$ $4$-chain in $\Delta_{
\beta +e_1 -e_2}$ with $\partial \tilde{\eta}= \partial \eta$.
\end{corollary}

{\it Proof.}
To prove the statement,  is sufficient to prove it when 
 $\eta$ is a $UFO^2_{t,5}$ for $t=1,...,5 $,
since $\eta$ is a sum  of $UFO^2_{t,5}$ for $t=1,...,5 $.
    Thus our statement  follows from Lemmas \ref{ufisemplici},
  \ref{ufo3} and \ref{ufo2}.
\hfill \includegraphics[scale=0.9]{box}

\bigskip

{\it Proof of Prop. \ref{step2}.}
We will show that if $\gamma$ is a $3$-cycle  in $\Delta_{b} $ with 
$\gamma = \partial \eta$ with  $\eta $ 
 in $\Delta_{b}  \cup \Delta_{b -e_1 + e_2} \cup ...
 \cup \Delta_{b -ke_1 +k e_2}$  for some $k \leq  b_1$, 
 then 
we can construct $\eta'$ in  $\Delta_{b}  \cup ... \cup 
\Delta_{b -(k-1) e_1 +(k-1) e_2}$ s.t. $\partial \eta' = \gamma$
(this, by induction on $b_1$, implies obviously Prop. \ref{step2}):

 let $\nu$ be the 
sum of the simplexes of $\eta$ in 
$\Delta_{b -ke_1 +k e_2 }$
 and not in $\Delta_{b -k e_1 +(k-1) e_2}$;
 $\partial \nu $ is in $\Delta_{
b -k e_1 +(k-1) e_2}$; by Corollary 
\ref{fine} $\partial \nu  = \partial \nu' $ for some 
$\nu'$ in $\Delta_{b -(k-1) e_1 +(k-1) e_2}$
 let  $\eta'= \eta- \nu + \nu' $; $\eta'$ is in 
$\Delta_{b}  \cup ... \cup 
\Delta_{b -(k-1) e_1 +(k-1) e_2}$
 and $\partial \eta' = \partial \eta= \gamma$.
\hfill \includegraphics[scale=0.9]{box}


{\small

\smallskip

{\bf Address: Elena Rubei, Dipartimento di Matematica ``U. Dini'',
via Morgagni 67/A, 
50134 Firenze, Italia.}
{\bf E-mail address: rubei@math.unifi.it}

\smallskip

{\bf 2000 Mathematical Subject Classification:} 
14M25, 13D02.}

\end{document}